\documentclass[preprint,12pt]{elsarticle}
\usepackage{amsmath}
\usepackage{amsfonts}
\usepackage{amssymb}
\usepackage{latexsym}
\usepackage{graphicx}
\usepackage{extarrows}
\usepackage{color}
\newtheorem{theorem}{\bf Theorem}[section]
\newtheorem{lemma}[theorem]{\bf Lemma}

\newtheorem{coro}[theorem]{\bf Corollary}

\newcommand{\R}{{\rm I\!R}}
\newcommand{\C}{\mathbb{C}}

\begin{document}

\title{Hardy-Hodge Decomposition of Vector Fields in $\mathbb{R}^n$ }
%\author{Laurent Baratchart \thanks{INRIA Sophia-Antipolis, France. Email:Laurent.Baratchart@inria.fr}, \ Pei Dang \thanks{Faculty of Information Technology, Macau
%University of Science and Technology, Macao, China. Email: pdang@must.edu.mo. Telephone: +853 88972823. The work was supported by Macao Science and Technology Development Fund, MSAR. Ref. 045/2015/A2}, \ Tao~Qian \thanks{Department of Mathematics, University of Macau, Macao (Via Hong
%Kong). Email: fsttq@umac.mo. Telephone: +853 83978547. Fax: +853 28838314. The work was supported by
%Multi-Year Research Grant (MYRG) MYRG116(Y1-L3)-FST13-QT, Macao Government FDCT 098/2012}}

\author[a]{Laurent Baratchart}

\author[b]{Pei Dang}

\author[c]{Tao Qian \footnote{Corresponding author. Tel: (+00853) +853 83978547. Fax: +853 28838314. Email: fsttq@umac.mo.}}

 \address[a]{INRIA, 2004 route de lucioles, 06902 Sophia-Antipolis Cedex, France. Email: Laurent.Baratchart@inria.fr}

\address[b]{Faculty of Information Technology, Macau
University of Science and Technology, Macao, China. Email: pdang@must.edu.mo. }

\address[c]{Department of Mathematics, University of Macau, Macao (Via Hong
Kong). Email: fsttq@umac.mo.}

\begin{abstract}
\textcolor{black}{We prove that a $\R^{n+1}$-valued vector field on $\R^n$ is the sum of
the traces of two harmonic gradients, one in each component
of $\R^{n+1}\setminus\R^n$, and of a $\R^n$-valued divergence free vector
field. We apply this to the description of vanishing potentials in divergence
form. The results are stated in terms of Clifford Hardy spaces, the
structure of which is important for our study.}
% Hardy spaces decomposition of scalar-valued functions in the Lebesgue $L^p$ spaces on a manifold has both the theoretical and applicable significance. In this paper we study similar decompositions but for para-victor fields defined on the Euclidean space $\mathbb{R}^n.$ Precisely, we show that a para-vector-field $f=f_0{\bf e}_0+f_1{\bf e}_1+\cdots +f_n{\bf e}_n$ in $L^p(\mathbb{R}^n), \ 1<p<\infty,$ where $f_k, k=0,1,...,n,$ are real-valued functions, can be decomposed as $f=f^++f^-+f^0,$ where $f^\pm$ are the non-tangential boundary limits on $\mathbb{R}^n$ of some two para-vector-valued  functions in, respectively, the para-vector valued Hardy spaces $H^p(\mathbb{R}^{n+1}_\pm),$ and $f^0$ is a divergence free vector-field on $\mathbb{R}^n.$
\end{abstract}

\maketitle

\section{Introduction}

\textcolor{black}{Decomposing a complex function on the line as the sum of the traces of
two holomorphic functions, one in each half plane cut out by the line, is
a classical topic from complex analysis that
lies at the root of many
developments in harmonic analysis.
Indeed, such a decomposition features the Hilbert transformation which is the
prototype of a Calder\`on-Zygmund operator, whose $C^{1,\alpha}$ and
$L^p$ boundedness
was historically the starting point of  elliptic regularity theory
\cite{Ga,Muskhelishvili,St,Torchinsky}.  This decomposition is also a
cornerstone
of solutions to Riemann-Hilbert problems, which are especially meaningful
in spectral theory  \cite{Fokas} and provided in recent
years striking advances in the theory of orthogonal polynomials \cite{Deift}.
Moreover,
it is instrumental for defining and
studying Hankel and Toeplitz operators, which
play a fundamental role in complex approximation and were successfully applied
to issues of basic importance in control and signal analysis
\cite{Nikolskii, Peller, Parfenov, Co}. Besides, in a Hilbertain context,
the decomposition was used to}
obtain sparse representations \textcolor{black}{of} analytic signals of scalar-valued signals in various classical contexts \textcolor{black}{(\cite{QWa} 
and \textcolor{black}{subsequent} papers \textcolor{black}{by these authors)}}.

Specifically, given a complex-valued \textcolor{black}{function}
$f\in L^p(\mathbb{R}), \ 1<p<\infty,$ one \textcolor{black}{has}
\begin{equation}
\label{decomp}
 f=f^++f^-,
\end{equation}
where
\[ f^\pm (x)=\lim_{y\to 0\pm}\frac{\pm 1}{2\pi i}\int_{-\infty}^\infty \frac{f(t)}{t-(x+iy)}dt,\]
where
$f^\pm$ are, respectively, non-tangential boundary limit functions of holomorphic functions of one complex variable in, respectively, the Hardy spaces $H^p(\mathbb{C}^\pm)$ \textcolor{black}{of the upper and lower half planes}.
The Hardy space functions are given by
\[ f^\pm (z)=\frac{\pm 1}{2\pi i}\int_{-\infty}^\infty \frac{f(t)}{t-z}dt,\quad z=x+iy, \ \pm y>0\textcolor{black}{,}\]
see \cite{Ga}.

\textcolor{black}{Now, under the standard identification $\C\sim\R^2$,
a holomorphic function may be regarded as the gradient of a harmonic function,
and this way \eqref{decomp} says that a $\R^2$-valued vector field on $\R$
of $L^p$ class
may be decomposed as the sum of the traces of two harmonic gradients,
arising from harmonic functions in the upper and lower half plane respectively.
The question that we raise in this paper is whether such a decomposition is
possible in higher dimension, namely whether
a vector field in $L^p(\R^n,\R^{n+1})$ is the sum of the traces of two
gradients of functions harmonic in the two half-spaces cut out by
$\R^n$ in $\R^{n+1}$.
The answer is no in general, but the next
best thing is that a decomposition  becomes possible if a third summand is
allowed, which takes the form of
a  divergence free vector field tangent to $\R^n$. This fact was observed in
\cite{BHLSW} when $n=2$, and used to characterize silent magnetization
distributions on a plane. We presently
carry this decomposition  over to every $n$. When projected onto $\R^n$,
it yields back the classical Hodge decomposition of a $L^p$ tangent
vector  field on $\R^n$ as the sum of a gradient and of a divergence free
component. This why we call our  decomposition of  $L^p(\R^n,\R^{n+1})$
vector fields the Hardy-Hodge decomposition.
}

\textcolor{black}{Formally  the decomposition can be surmised
from Hodge theory for
1-currents supported on a hypersurface in ambient space 
\cite[Sec. 2.8]{Harris}, 
but the estimates needed to control
$L^p$-norms of the objects involved pertain to the Calder\`on-Zygmund theory.}
\textcolor{black}{In this connection, it would be pedantic to
introduce currents to speak of vector fields on linear submanifolds, but it
is convenient to use the formalism of Clifford analysis, which provides us
with a substitute for complex variables and is well adapted to 
handle higher dimensional singular integrals}. In fact, Clifford analysis is 
also suited to extend the
result to vector fields on more general sub manifolds, 
although 
such a generalization lies beyond the scope of the present paper.

\textcolor{black}{The latter is organized as follows. In Section \ref{prelim}
we recall some basic facts from Clifford analysis and Clifford
Hardy spaces, most of which can be found in \cite{GM}, and we study the
structure of boundary function in detail, along with density properties
of rational-like functions. In Sections \ref{sHH} and \ref{var}, we
prove the Hardy-Hodge decomposition and some variants thereof. Finally, in
Section  \ref{sis}, we discuss an application to non-uniqueness
for inverse potential problems in divergence form.}

% Such decomposition of Functions on manifolds are mathematically well known but have recently found significant applications in signal analysis in relation to the so called analytic signals \cite{Co}. In particular,  Hardy spaces are reproducing kernel Hilbert spaces. One, therefore, can obtain sparse representations of Hardy space functions, and, subsequently, sparse representations of scalar-valued functions in various classical contexts (See \cite{QWa} and Qian et al recent publications).

\section{Preliminaries}
\label{prelim}
\textcolor{black}{Let $n\geq 3$ be an integer and $\Phi$ be either the real field $\R$ or the complex field $\C$.
Hereafter,
we put $L^p(\R^n,E)$ for the familiar Lebesgue space of functions on
$\R^n$ with values in a Banach space $E$ (typically $E=\Phi^m$)
whose norm to the $p$-th power is integrable, and
we often write $L^p(\R^n)$ for
simplicity if $E$ is understood from the context.}

We adopt standard notations in \textcolor{black}{ Clifford analysis, see \cite{GM}}.
In particular, we put \textcolor{black}{$Cl(n,{\Phi})$}  to
denote the Clifford algebra generated over
\textcolor{black}{$\Phi$} by ${\bf e}_1, ..., {\bf e}_n$
\textcolor{black}{ with} ${\bf e}_0=1$ and
${\bf e}_i{\bf e}_j+ {\bf e}_j{\bf e}_i=-2\delta_{ij}, \ i,j=1,...,n$.
\textcolor{black}{We indicate with $\mathcal{S}$ the collection of subsets of
$\{1,\cdots,n\}$. Then,
the elements of the canonical basis of
 $Cl(n,\Phi)$, viewed as a vector space
over $\Phi$, are denoted as ${\bf e}_S$, $S\in\mathcal{S}$, where
$e_S={\bf e}_{j_1}\cdots {\bf e}_{j_k}$ if}
$S=\{1\leq j_1<\cdots <j_k\leq n\}$.   \textcolor{black}{A generic member of
$Cl(n,{\Phi})$ can thus be written as $x=\Sigma_{S\in\mathcal{S}}\,x_S {\bf e}_S$ with
$x_S\in\Phi$.}
When $S$ is empty,
we \textcolor{black}{write} ${\bf e}_\emptyset={\bf e}_0=1$. \textcolor{black}{The conjugate of $x$, denoted as $\overline{x}$, is defined to be
$\Sigma_{S\in\mathcal{S}}\,(-1)^{|S|}x_S {\bf e}_S$, where $|S|$ indicates the cardinality of $S$.
By convention,} a $0$-form is a scalar. \textcolor{black}{A $k$-form is a sum
$\Sigma_{S\in\mathcal{S}_k}\,x_S {\bf e}_S$ where $\mathcal{S}_k$ indicates
those members of $\mathcal{S}$ with cardinality $k$.}
\textcolor{black}{Clearly, $Cl(n,\Phi)$} is a $2^n$-dimensional linear space
\textcolor{black}{over $\Phi$}. A $1$-form \textcolor{black}{is also}
called a \emph{vector}, \textcolor{black}{denoted  with underscore:} $\underline{x}=x_1{\bf e}_1+\cdots +x_n{\bf e}_n$. \textcolor{black}{Clifford vectors are identified with Euclidean vectors in $\R^n$. The sum of
a $0$-form and  a $1$ form} is called a \emph{para-vector}, \textcolor{black}{and if $x$ is a paravector we let $\underline{x}$ be its
vector part:} $x=x_0{\bf e}_0+\underline{x}$. \textcolor{black}{This is consistent with our previous notation for vectors}. The norm of $x\in Cl(n,\Phi)$ is defined to be $|x|=(\sum_{\textcolor{black}{S\in\mathcal{S}}} |x_S|^2)^{1/2}$,
\textcolor{black}{which derives from} the inner product
$\langle x,y\rangle=\sum_{\textcolor{black}{S\in\mathcal{S}}} x_S\overline{y}_S.$ If both $x,y$
are para-vectors, then \textcolor{black}{their} Clifford product $xy=-\langle x,y\rangle +x\wedge y,$ where
%for $n=3,$ for instance,
the exterior product $x\wedge y$ \textcolor{black}{is a 2-form similar to
the exterior  product of differential forms from geometry:
\[
x\wedge y=\Sigma_{j<k}(x_jy_k-x_ky_j){\bf e}_j{\bf e}_k,\qquad
x=\Sigma_{S\in\mathcal{S}}\,x_S {\bf e}_S,\quad y=\Sigma_{S\in\mathcal{S}}\,y_S {\bf e}_S.\]}
%$+(x_2y_3-x_3y_2){\bf e}_2{\bf e}_3+(x_3y_1-x_1y_3){\bf e}_3{\bf e}_1.$
\textcolor{black}{Let} $\mathbb{R}^{n}_{1} =\{x=x_0+\underline{x}\ |\  x_0\in \mathbb{R}, \ \underline{x}\in \mathbb{R}^n\},$ $\mathbb{R}^{n}_{1,\pm} =\{x=x_0+\underline{x}\ |\ \pm x_0>0, \ \underline{x}\in \mathbb{R}^n\},$ and
$\mathbb{R}^{n+1}_\pm =\{ x=\underline{x}+x_{n+1}{\bf e}_{n+1}\ |\ \pm x_{n+1}>0\}$, \textcolor{black}{noting that $\mathbb{R}^{n}_{1}, \mathbb{R}^{n}_{1,\pm}$ are included in $Cl(n,\R)$ while $\mathbb{R}^{n+1}_\pm $ is contained in
$\R^{n+1}\subset Cl(n+1,\R)$ (via the identification of Clifford vectors with Euclidean vectors). For $k=0,...,n+1$, we introduce the partial derivatives}  $\partial_k=\partial/\partial_{x_k}$ \textcolor{black}{and subsequently  we
define}
\begin{equation}
\label{defop}
\begin{array}{lcl}
&D_0=\partial_0,\quad  &D_n={\bf e}_1\partial_1+\cdots {\bf e}_n\partial_n,\\
&D=D_0+D_{n}, \quad  &D_{n+1}={\bf e}_1\partial_1+\cdots {\bf e}_n\partial_n +{\bf e}_{n+1}\partial_{n+1}.
\end{array}
\end{equation}
\textcolor{black}{A $Cl(n,\Phi)$ or $Cl(n+1,\Phi)$-valued function $f$ such that
$Df=0\ (fD=0)\ $ or $D_{n+1}f=0$ ($fD_{n+1}=0$) on an open set of
$\Omega\subset\R^{n+1}$} is called \emph{left-monogenic} (\emph{right-monogenic}) \textcolor{black}{on $\Omega$. By convention,
coordinates in the  case of $D$ are denoted by
$x_0,\cdots,x_n$, whereas in the case of $D_{n+1}$ they
are written $x_1,\cdots,x_{n+1}$}.  If a function is both left- and right-monogenic, we call \textcolor{black}{ it}  \emph{two-sided-monogenic}.
\textcolor{black}{Let us stress that, when applying the differential operators
\eqref{defop}, the partials $\partial_j$ commute with the $e_k$ but the $e_j$
do not, so that it generally
matters whether the operator gets applied from the left or
the right.}

\textcolor{black}{Note that $(D_0-D_n)D=\Delta$ (resp. $D_{n+1}^2=-\Delta$)
where $\Delta=\Sigma_{j=0}^n\partial^2_{x_j}$
(resp.  $\Delta=\Sigma_{j=1}^n\partial^2_{x_j}$ ) is the ordinary Laplacian.
Therefore left or right monogenic functions have harmonic components,
in particular they are real analytic on $\Omega$ and there is no difference
being monogenic in the distributional or in the strong sense.
When $f$ is para-vector valued in $Cl(n,\Phi)$ and we write  $f=f_0+f_1{\bf e}_1+\cdots+f_n {\bf e}_n$, it is readily checked that
$Df=0$ if and only if
\begin{equation}
\label{paramon}
\partial_0 f_0=\Sigma_{j=1}^n\partial_j f_j\quad \mbox{\rm with}\quad\partial_jf_k=\partial_k f_j\quad\mbox{\rm and}\quad
\partial_0f_j=-\partial_j f_0\quad \mbox{\rm for}\quad
 1\leq j<k\leq n,
\end{equation}
and similarly when $fD=0$. In particular, a} vector or para-vector valued
\textcolor{black}{function which is left-monogenic must  be also right-monogenic.
In the same manner, when $f$ is vector valued in $Cl(n+1,\Phi)$ with
$f=f_1{\bf e}_1+\cdots+f_{n+1} {\bf e}_{n+1}$, we have that
$D_{n+1}f=0$ if and only if
\begin{equation}
\label{paramon1}
\Sigma_{j=1}^{n+1}\partial_j f_j=0\quad \mbox{\rm and}\quad\partial_jf_k-\partial_k f_j=0\quad \mbox{\rm for}\quad
 1\leq j<k\leq {n+1},
\end{equation}
and the same if $fD_{n+1}=0$. Thus, vector valued left-monogenic functions
are  right-monogenic}.

\textcolor{black}{By \eqref{paramon}, a para-vector valued function
$f=f_0+f_1{\bf e}_1+\cdots+f_n {\bf e}_n$, where the $f_j$ are real valued, is
monogenic if and only if $(-f_0,f_1,\cdots,f_n)$ is a harmonic gradient,
meaning that it is the gradient of a harmonic function.
The components of a harmonic gradient are sometimes referred to as
a conjugate harmonic system, or a Riesz system of functions, {\it cf.}
\cite{SW}. When
$\R_1^n$ gets identified with $\R^{n+1}$, the fact that $(-f_0,f_1,\cdots,f_n)$
is a harmonic gradient  amounts to saying that
$\overline{f}=f_0-f_1{\bf e}_1-\cdots-f_n {\bf e}_n$ is a harmonic gradient.
Likewise, it follows from \eqref{paramon1} that
a vector valued function in $Cl(n+1,\R)$, say
$f=f_1{\bf e}_1+\cdots+f_{n+1} {\bf e}_{n+1}$ is monogenic if and only if
$(f_1,\cdots,f_{n+1})$ is a harmonic gradient. Identifying
vectors in $Cl(n+1,\R)$ with $\R^{n+1}$, we simply say in this case
that $f$ is a harmonic gradient.
}

Let $g: \mathbb{R}^{n}_{1,\pm}\to \mathbb{R}^{n}_{1}$. For $1<p<\infty,$ we say $g$ belongs to the Hardy space $H^p(\mathbb{R}^{n}_{1,\pm}, \mathbb{R}^n_1)$ if $Dg=0$ in \textcolor{black}{$\mathbb{R}^{n}_{1,\pm}$} and
\begin{eqnarray}\label{firste}  \|g\|^{\textcolor{black}{p}}_{H_{\textcolor{black}{\pm}}^p}\textcolor{black}{\triangleq}\sup_{\pm x_0>0} \int_{\mathbb{R}^n}|g(x_0+\underline{x})|^p d\underline{x}<\infty.\end{eqnarray}
We \textcolor{black}{refer to} the \textcolor{black}{above Hardy} spaces as
\textcolor{black}{being of  para-vector  \textcolor{black}{type},
or \textcolor{black}{also of} \emph{inhomogeneous} type}.
\textcolor{black}{Thus, $g\in H^p(\mathbb{R}^{n}_{1,\pm}, \mathbb{R}^n_1)$ if and only if $\overline{g}$ is a harmonic gradient
which moreover satisfies} the $p$-norm boundedness (\ref{firste}) in the relevant half space.
\textcolor{black}{Equivalently, since each of the functions composing a conjugate harmonic
system is
harmonic, it follows from (\ref{firste}) and standard estimates
on harmonic functions (see {\it e.g.}
\cite[Ch. II, Thm. 3.7 \&  Eqn. (3.18)]{SW})
that a para-vector valued monogenic function in $\R_{1,\pm}^n$ lies in
$H^p(\mathbb{R}^{n}_{1,\pm}, \mathbb{R}^n_1)$ if and only if
the
nontangential
maximal function given by
\begin{equation}
\label{defntm}
\mathcal{M}_\alpha g(\underline{x})\triangleq\sup_{x_0+\underline{x}\in
\Gamma_\alpha(\underline{x})}|g(x_0+\underline{x})|
\end{equation}
lies in $L^p(\R^n,\R)$ with equivalence of norms: $\|g\|_{H^p_{\pm}}\leq
\|\mathcal{M}_{\alpha} g\|_{L^p(\R^n)}\leq C_\alpha\|g\|_{H^p_{\pm}}$. Here,
to each $\alpha>0$ and $\underline{x}\in\R^n$, the notation
$\Gamma_\alpha(\underline{x})$ stands for the cone
$$\Gamma_\alpha(\underline{x})=\{y_0+\underline{y}\in \R_{1,\pm}^n, \ |\underline{y}-\underline{x}|< \alpha|y_0|\},$$
and the precise value of $\alpha$ is irrelevant except that the constants will depend on it. }
%of this gradient is bounded in
%$L^p(\mathbb{R}^n)$} (\cite{SW}, \cite{St}, \cite{LMcQ}).

\textcolor{black}{Likewise, for} $g: \mathbb{R}^{n+1}_{\pm}\to \mathbb{R}^{n+1}\subset Cl(n+1,\R)$ \textcolor{black}{and}  $1<p<\infty,$ we say \textcolor{black}{that} $g$ belongs to the Hardy space $H^p(\mathbb{R}^{n+1}_\pm, \mathbb{R}^{n+1})$ if $D_{n+1}g=0$ in the
half \textcolor{black}{space $\R_\pm^{n+1}$}  and
\begin{eqnarray}\label{gradiants}  \textcolor{black}{\|g\|_{H^p_{\pm,h}}}=\sup_{\pm x_{n+1}>0} \int_{\mathbb{R}^n}|g(\underline{x}+x_{n+1}{\bf e}_{n+1})|^p d\underline{x}<\infty.\end{eqnarray}
Note that, in the previous equation, $\underline{x}$ refers to a vector in $Cl(n,\R)$ viewed as a vector in $Cl(n+1,\R)$ whose $(n+1)$-st component is zero.
We refer \textcolor{black}{to the latter} Hardy spaces as \textcolor{black}{being of} vector \textcolor{black}{type} or homogeneous  type, \textcolor{black}{which is the reason for the subscript ``$h$'' in the notation for the norm. Thus, we have
that $f=f_1{\bf e}_1+\cdots+f_{n+1} {\bf e}_{n+1}$ lies in  $H^p(\mathbb{R}^{n+1}_\pm, \mathbb{R}^{n+1})$ if and only if it
is a harmonic gradient in $\mathbb{R}^{n+1}_\pm$ satisfying the $p$-boundedness condition \eqref{gradiants}. The
latter is again equivalent to the $L^p$ boundedness on $\mathbb{R}^n$ of the nontangential maximal function
\begin{equation}
\label{defntmh}
\mathcal{M}_{\alpha,h} g(\underline{x})\triangleq\sup_{\underline{x}+x_{n+1}{\bf e}_{n+1}\in
\Gamma_{\alpha,h}(\underline{x})}|g(\underline{x}+x_{n+1}{\bf e}_{n+1})|
\end{equation}
where this time
$$\Gamma_{\alpha,h}(\underline{x})=\{\underline{y}+y_{n+1}{\bf e}_{n+1}\in
\R_{\pm}^{n+1}, \ |\underline{y}-\underline{x}|< \alpha|y_{n+1}|\}.$$
In fact the passage from non homogeneous to homogeneous Hardy spaces is
rather mechanical, trading $x_0$ and $1$ for $x_{n+1}$ and ${\bf e}_{n+1}$
while changing $g_0$ into $-g_{n+1}$.}

\textcolor{black}{Next, recall the local Fatou theorem asserting that a harmonic function
in $\R_+^{n+1}$ which is non-tangentially bounded at almost every point of
a set $G\subset\R^n$ has a non-tangential limit at almost every point of $G$;
here, nontangential refers to the fact that bounds and limits are seeked in
cones $\Gamma_{\alpha,h}(\underline{x})$ for arbitrary but fixed $\alpha>0$,
see \cite[Thm. 3.19]{SW}. In view of \eqref{firste} and \eqref{gradiants},
it follows from the local Fatou theorem
 that each component of a Hardy function (resp. homogeneous Hardy function) has a nontangential limit at almost every point of
$\R^{n}\subset Cl(n,\R)$ (resp. $\R^{n}\times\{0\}\subset\R^{n+1}\subset Cl(n+1,\R)$). This defines boundary values for such functions. Now,
it is an important and peculiar property of
left or right monogenic functions that they can be recovered as Cauchy integrals of their boundary values, see  \cite[Cor. 3.20]{GM} and \cite{KQ, LMcQ}. Specifically,
% an important property of the Clifford Hardy spaces that each $H^p(\mathbb{R}^{n}_{1,\pm})$ function $g$ has non-tangential boundary limits almost everywhere on $\mathbb{R}^n,$  still denoted by $g,$ and the $g\in H^p(\mathbb{R}^{n}_{1,\pm})$ may be expressed by the Cauchy integral of its boundary limit function $g\in L^p(\mathbb{R}^n)$
let us discuss the case of left monogenic $Cl(n,\R)$-valued functions on
$\R_{1,+}^n$  as a prototypical example. If $g$ is such a function
and if $\mathcal{M}_\alpha g\in L^p(\R^n)$ for some $p\in(1,\infty)$, then
$g$ has nontangential limit a.e. on $\R^n$
 since its components are harmonic functions with
$L^p$ nontangential maximal function, and this nontangential limit clearly lies in $L^p(\R^n)$ because it is dominated by the nontangential maximal function.
Then, denoting
the nontangential limit by $g$ again \textcolor{black}{and letting $\omega_n$ designate the volume of
the unit sphere of dimension $n$,},  we get
 that $g=C^+ g$, where}

\begin{eqnarray}
\label{defCauchyp}
C^+g(x)&=&\frac{1}{\omega_n}\int_{\mathbf{R}^n}\frac{\overline{\underline{y}-x}}
{|\underline{y}-x|^{n+1}}(-e_0)g(\underline{y})d\underline{y}\\
\nonumber
&=&\frac{1}{\omega_n}\int_{\mathbf{R}^n}
\frac{x_0}
{|x-\underline{y}|^{n+1}}g(\underline{y})d\underline{y}+\frac{1}{\omega_n}\int_{\mathbf{R}^{n}}
\frac{x_1-y_1}
{|x-\underline{y}|^{n+1}}(-\mathbf{e}_1)g(\underline{y})d\underline{y}+\cdots + \\
\nonumber
&&+\frac{1}{\omega_n}\int_{\mathbf{R}^{n}}
\frac{x_n-y_n}
{|x-\underline{y}|^{n+1}}(-\mathbf{e}_n)g(\underline{y})d\underline{y}, \ \ \ \ \ \ \ \ x=x_0\mathbf{e}_0+x_1\mathbf{e}_1+\cdots +x_n\mathbf{e}_n\in \mathbf{R}^n_{1,+},
\end{eqnarray}
\textcolor{black}{with} $\omega_n$ the surface measure of the $n$-dimensional unit sphere. \textcolor{black}{Here,
the presence of $-e_0$  in the
definition of $C^+$ is because  $-e_0$ is the exterior unit normal to
$\R^n_{1,+}$, see \cite{GM,LMcQ}.}

\textcolor{black}{Conversely, if $g$ is any $Cl(n,\R)$-valued
function in $L^p(\R^n)$, the formula for $C^+ g$ makes good sense and defines a left monogenic function in $\R_{1,+}^n$ since the kernel is left monogenic (as is easily verified).} Invoking the Plemelj formula \cite{LMcQ},
\textcolor{black}{we get at}  almost all points \textcolor{black}{$\underline{x}\in\R^n$ that} there exist non-tangential \textcolor{black}{limits  of $C^+g$} denoted
\textcolor{black}{with a} curly $\mathcal{C}^+g,$ \textcolor{black}{which are given by}
\begin{eqnarray}\label{Cauchy PI}
\mathcal{C}^+g(\underline{x})=\frac{1}{2}[g(\underline{x})+Hg(\underline{x})]=\frac{1}{2}(I+H)g(\underline{x}),
\end{eqnarray}
where
\begin{eqnarray*}
Hg(\underline{x})&=&\frac{\textcolor{black}{2}}{\omega_n}\lim_{\varepsilon\to 0^+}\int_{|\underline{x}-\underline{y}|>\varepsilon}
\frac{x_1-y_1}
{|\underline{x}-\underline{y}|^{n+1}}(-\mathbf{e}_1)g(\underline{y})d\underline{y}
+\cdots + \frac{\textcolor{black}{2}}{\omega_n}\lim_{\varepsilon\to 0^+}\int_{|\underline{x}-\underline{y}|>\varepsilon}
\frac{x_n-y_n}
{|\underline{x}-\underline{y}|^{n+1}}(-\mathbf{e}_n)g(\underline{y})
d{\underline y}\\
&\triangleq& \sum_{k=1}^n(-\mathbf{e}_k)R_k(g)(\underline{x}),
\end{eqnarray*}
\textcolor{black}{where}
$$R_k(g)(\underline{x})=\frac{\textcolor{black}{2}}{\omega_n}\lim_{\varepsilon\to 0^+}\int_{|\underline{x}-\underline{y}|>\varepsilon}\frac{x_k-y_k}
{|\underline{x}-\underline{y}|^{n+1}}g(\underline{y})d\underline{y}$$
is the $k$-th Riesz \textcolor{black}{transformation} of $g, \ k=1,...,n$.
\textcolor{black}{As an operator, $R_k$ has multiplier
$-i\kappa_k/|\kappa|$ in the Fourier domain ($\kappa=(\kappa_1,
\cdots,\kappa_n)$ denoting the
Fourier variable)
and it maps $L^p(\R^n,\R)$ into
itself for $1<p<\infty$, see \cite[Ch. III, Secs. 4.3-4.4]{St}.}
\textcolor{black}{The operator $H$ defines} the Hilbert transformation \textcolor{black}{in the present} context,
\textcolor{black}{and it satisfies} the relation $H^2=I,$ where $I$ is the identity.

We note that \textcolor{black}{our}  definition \textcolor{black}{of} $H$ is
consistent with the \textcolor{black}{classical definition of the
Hilbert transform in dimension $1$, }
given \textcolor{black}{{\it e.g.} in}
 \cite{Be}.  \textcolor{black}{For if we designate the latter with $H_1$},
\textcolor{black}{we get} when $n=1,$ with ${\bf e}_1=-i,$ \textcolor{black}{that}
 $$Hg=-{\bf e}_1 \frac{1}{\pi}\left({\rm p.v.}\frac{1}{(\cdot )}\ast\right) g=i \frac{1}{\pi}\left({\rm p.v.}\frac{1}{(\cdot )}\ast\right) g\triangleq iH_1g,$$
\textcolor{black}{where ``p.v.'' indicates the principal value.}
\textcolor{black}{Thus, in view of} the well known relation $H_1^2=-I,$ we
\textcolor{black}{have that} $H^2=(iH_1)^2=I,$ \textcolor{black}{which is}
consistent with the \textcolor{black}{$n$-dimensional}  case.
\textcolor{black}{Similar considerations apply to functions in
$H^p(\mathbb{R}^{n}_{1,-}, \mathbb{R}^n_1)$, only trading $C^+$ for
its opposite (because the outer normal to
$\R^n_{1,-}$ is ${\bf e}_0$ and not $-{\bf e}_0$ ), and letting $x$ now
range  over $\mathbb{R}^{n}_{1,-}$. Letting this time $\mathcal{C}^-g$
stand for the nontangential limit on $\R^n$, this results in the
Plemelj formula:
\begin{eqnarray}\label{Cauchy PIm}
\mathcal{C}^-g(\underline{x})=
%\frac{1}{2}[g(\underline{x})-Hg(\underline{x})]=
\frac{1}{2}(I-H)g(\underline{x}).
\end{eqnarray}
}

In the sequel we denote by ${\rm Sc}\{x\}$ \textcolor{black}{the} scalar part of a Clifford number $x\in Cl(n,\Phi),$ \textcolor{black}{which is}
the $0$-form of $x,$ and \textcolor{black}{by} ${\rm Nsc}\{x\}$ \textcolor{black}{the} non-scalar part of $x,$ \textcolor{black}{which is} the sum of all the $k$-forms of $x,\ k=1,...,2^n.$

%We call a (para-) vector-valued function $g$ \textcolor{black}{the} harmonic gradient of some scalar-valued harmonic function $h$ if
%$g=\triangledown h$\textcolor{black}{; here, in the para-vector valued case,
%$\R_1,n$ gets naturally identified with $\R^{n+1}$, as usual.
\textcolor{black}{We pointed out already that each function in
$H^p(\mathbb{R}^{n}_{1,\pm}, \mathbb{R}^n_1)$ is naturally associated to
the conjugate of a harmonic gradient.
In the next lemma,
we identify the latter as being the gradient of a Newton potential, and we
describe the boundary values of
Hardy functions.}

\begin{lemma}
\label{bvH}
If $g=g_0{\bf e}_0+\cdots +g_n{\bf e}_n\in H^p(\mathbb{R}^{n}_{1,\pm}, \mathbb{R}^n_1),$ where each $g_k$ is real valued, then its non-tangential
\textcolor{black}{limit on $\R^n$}, still denoted as $g,$ satisfies $g=(I\pm H)g_0,$ that is, $g_k=\textcolor{black}{\mp}R_kg_0$ \textcolor{black}{for} $ k=1,...,n.$ \textcolor{black}{Conversely, each function on $\R^n$
of the form $(I\pm H) \varphi$, with $\varphi\in L^p(\R^n,\R)$, is the
nontangential limit of a function in $H^p(\mathbb{R}^{n}_{1,\pm}, \mathbb{R}^n_1)$ and the Hardy norm is equivalent to $\|\varphi\|_{L^p(\R^n)}$.}
Moreover, the conjugate of each \textcolor{black}{$g\in H^p(\mathbb{R}^{n}_{1,\pm},\R_1^n)$}
%or $H^p(\mathbb{R}^{n+1}_{\pm})$
is \textcolor{black}{the} harmonic gradient \textcolor{black}{ of the Newton potential of $\mp2/(n-1)$ times its scalar part $g_0$, namely:}
\textcolor{black}{
\begin{equation}
\label{cNewton}
g(x)=\partial_0\mathcal{N}_{g_0}(x)-\Sigma_{k=1}^{n}\partial_k\mathcal{N}_{g_0}(x){\bf e}_k,
\qquad \mathcal{N}_{g_0}(x)\triangleq
\mp\frac{2}{(n-1)\omega_n}\int_{\R^n}\frac{g_0(\underline{y})}{|x-\underline{y}|^{n-1}}d\underline{y},\qquad x\in\R_{1,\pm}^n.
\end{equation}}
\end{lemma}

\noindent{\bf Proof.}  \textcolor{black}{By the Cauchy-Clifford formula,
the function $g$ in the upper-half space may be expressed as the Cauchy integral over $\mathbb{R}^n$ of its boundary limit function. Therefore the}
Plemelj theorem implies \textcolor{black}{that}
$(1/2)(I\pm H)g=g$\textcolor{black}{, hence $\pm Hg=g$.}
By comparing the terms \textcolor{black}{of degree $1$}
on the left- and the right-hands of the above identity we \textcolor{black}{get
$\mp{\bf e}_kR_kg_0={\bf e}_kg_k$ for  $ k=1,...,n$, ensuing that} $g=(I\pm H)g_0$.  \textcolor{black}{Invoking the Plemelj theorem again, this implies that
the Cauchy integral of $2g_0$ is a left monogenic function whose nontangential
boundary value is $g$. Hence, by the Cauchy formula, this function must coincide with $g$ on $\R_{1,\pm}^n$.
Now, when identifying $\R_{1}^n$ with
$\R^{n+1}$, the conjugate of the Cauchy kernel is
$-1/(n-1)$ times
the gradient} of the Newton \textcolor{black}{kernel $1/|x-.|^{n-1}$, and taking gradient
commutes} with the integration, as the integrand does not have singularity,
\textcolor{black}{so} we conclude that \textcolor{black}{the conjugate
of} $g$ is \textcolor{black}{the gradient of the Newton potential of $\mp2g_0/(n-1)$ (the $\mp$ arises because of the presence of $\mp{\bf e}_0$ in the Cauchy integral).
Conversely,  by \eqref{Cauchy PI}, a function  on $\R^n$
of the form $(I\pm H)\varphi$, where $\varphi$ is scalar valued
in $L^p(\R^n)$, is the nontangential limit of the
Cauchy integral $C^\pm (2\varphi)$ which is indeed
para-vector valued an monogenic. To see that it lies in
$H^p(\mathbb{R}^{n}_{1,\pm},\R_1^n)$, observe from what precedes that
its conjugate is the gradient of the Newton potential $\mathcal{N}_\varphi$
of $\mp2\varphi/(n-1)$.
In particular, by inspection of formula \eqref{cNewton}
(where $g_0$ is set to
$\varphi$), we find that
${\rm Sc}\{(C^\pm\varphi)(x_0,\underline{x})\}=(P_{x_0}*\varphi)(\underline{x})$
at every $(x_0,\underline{x})\in\R^n_{1,\pm}$, where the symbol ``$*$'' indicates convolution and
$P_{x_0}$ is the Poisson kernel at level $x_0$:
\begin{equation}
\label{Poissondef}
P_{x_0}(\underline{x})=\frac{2}{\omega_n}\frac{|x_0|}{\left(x_0^2+|\underline{x}|^2\right)^{(n+1)/2}},\qquad \underline{x}\in\R^n.
\end{equation}
Since $P_{x_0}$ has unit norm in $L^1(\R^n)$ for all $x_0$,
it follows that
\begin{equation}
\label{borP}
\left\|{\rm Sc}\left\{(C^\pm\varphi)(x_0,.)\right\}\right\|_{L^p(\R^n)}\leq\|\varphi\|_{L^p(\R^n)},\quad\forall x_0,
\end{equation}
implying that the scalar part of $C^\pm\varphi$ meets the
$p$-boundedness condition \eqref{firste}. To show that the vector part also
satisfies this condition, let us work on $\R_{1,+}^n$ as the
argument on  $\R_{1,-}^n$ is similar. Fix $z_0>0$ and consider the para-vector valued function on $\R_{1,+}^n$ given by
$F(x_0,\underline{x})= (C^+\varphi)(z_0+x_0,\underline{x})$. Clearly it is
monogenic, and we get upon applying H\"older's inequality to
\eqref{defCauchyp} that $\|F\|_{H^p_+}\leq c$ for some constant
$c=c(z_0)$, see definition \eqref{firste}. Hence
$F\in H^p(\mathbb{R}^{n}_{1,\pm},\R_1^n)$, and it is obvious that its
nontangential limit on $\R^n$ is $C^+\varphi(z_0,.)$. thus,
by the previous part of the proof, it holds that
\[(C^+\varphi)(x_0,.)=
(I+H){\rm Sc}\{F(0,.)\}=(I+H)
{\rm Sc}\left\{(C^+\varphi)(x_0,.)\right\}.
\]
In view of the definition of $H$ ({\it cf.} equation following
\eqref{Cauchy PI}) and the $L^p$-boundedness of Riesz transforms,
we now deduce from \eqref{borP} that
$C^+\varphi$ satisfies \eqref{firste}, as desired. We also proved
that $\|C^+\varphi\|_{H^p_+}\leq c\|\varphi\|_{L^p(\R^n)}$ for some constant $c$ independent of $\varphi$,
and since $\|\varphi\|_{L^p(\R^n)}$ is obviously
less than the $L^p$ norm of the non-tangential maximal function which itself
is equivalent to the Hardy norm, as pointed out after \eqref{defntm},
} the proof is complete.

\textcolor{black}{When dealing with functions in
$H^p(\mathbb{R}^{n+1}_\pm, \mathbb{R}^{n+1})$, the Cauchy formula is the
same except that $\mp{\bf e}_0$ gets replaced by $\mp{\bf e}_n$
(the outer normal to $\R_\pm^{n+1}$), and in the Plemelj formula
$H$ is changed into
$H{\bf e}_{n+1}=\Sigma_{k=1}^{n}(-{\bf e}_k{\bf e}_{n+1})R_k$, see \cite[Ch. 2, Sec. 5]{GM}.
An argument analogous to the previous one shows that
$f=f_1{\bf e}_1+\cdots+f_{n+1}{\bf e}_{n+1}$ lies in
$H^p(\mathbb{R}^{n+1}_\pm, \mathbb{R}^{n+1})$ if and only if
 $f_k=\textcolor{black}{\pm}R_kf_{n+1}$ \textcolor{black}{for} $ k=1,...,n.$ and
that
$f$ is $D_{n+1}\mathcal{N}$,
where $\mathcal{N}$ is
the Newton potential of $\mp2g_{n+1}/(n-1)$}.

\textcolor{black}{Observe that
Lemma \ref{bvH} and its analog for homogeneous spaces
entail that $H^p(\R^n_{1,\pm},\R_1^n)$ (resp.
$H^p(\R^{n+1}_\pm,\R^{n+1})$) is a
Banach space isomorphic to $L^p(\R^n,\R)$,
with norm equivalent to the $L^p$-norm of the trace
of the scalar part (resp. $n+1$-st component). Observe also from this, since
$L^q(\R^n,\R)\cap L^p(\R^n,\R)$ is dense in $L^p(\R^n,\R)$ for any
$q\in(1,\infty)$, that $H^p(\mathbb{R}^{n+1}_\pm, \mathbb{R}^{n+1})\cap H^q(\mathbb{R}^{n+1}_\pm, \mathbb{R}^{n+1})$ is dense in
$H^p(\mathbb{R}^{n+1}_\pm, \mathbb{R}^{n+1})$.
}

\textcolor{black}{Lemma \ref{bvH} easily implies
a result which is of interest in its own right and parallels
the density of rational functions in holomorphic Hardy spaces
of index $p\in(1,\infty)$ on the half-plane \cite{Ga}. Note that
rational functions with simple poles are conjugate of gradients of
logarithmic potentials of discrete  measures with finite support.
In the present context, analogs of rational functions with simple poles
are conjugates of gradients of Newton potentials of discrete
measures with finite support.
Specifically, if we let
\[
R_x(y)=\overline{\nabla_y\left(\frac{1}{\omega_n|x-y|^{n-1}}\right)}=
\frac{n-1}{\omega_n}\left(\frac{x_0-y_0}{|x-y|^{n+1}}-
\sum_{j=1}^n\frac{x_j-y_j}{|x-y|^{n+1}}{\bf e}_j\right),\quad
y\in\R^n_{1,+},\ x\in\R^n_{1,-},
\]
then $R_x\in H^p(\R^n_{1,+},\R^n_1)$ as a function of $y$ for fixed $x$,
and we have the following result.
\begin{coro}
The span of $\{R_x\}_{x\in\R^n_{1,-}}$ is dense in $H^p(\R^n_{1,+},\R^n_1)$
for $1<p<n$.
\end{coro}
\noindent{Proof.}
It follows from Lemma \ref{bvH} that $H^p(\R^n_{1,+},\R^n_1)$ is isomorphic
to $L^p(\R^n,\R)$ with equivalence of norms, the isomorphism being
\begin{equation}
\label{isoHL}
L^p(\R^n,\R)\ni h\ \mapsto \ C^+\Bigl(h-\Sigma_{k=1}^n (R_kh){\bf e}_k
\Bigr)\in H^p(\R^n_{1,+},\R^n_1).
\end{equation}
The inverse image of $R_x$ under this isomorphism is $(1-n)/2$ times
the Poisson kernel
$P_{x_0}(\underline{x}-\underline{y})$ defined in \eqref{Poissondef}.
Thus, by the Hahn-Banach theorem, the asserted density is equivalent to the
fact that no nonzero function in $L^{p'}(\R^n),\R)$, with $1/p+1/p'=1$,
can have vanishing Poisson integral. This, however, drops out
immediately from the property that the
Poisson kernel is an approximate identity\textcolor{black}{,
thereby achieving} the proof.}

\textcolor{black}{When} saying that a vector-valued function
\textcolor{black}{$f=f_1{\bf e}_1+\cdots+f_n{\bf e}_n$}
on ${\mathbb{R}^n}$ is divergence free, we mean \textcolor{black}{that}
${\rm div}{f}=\sum_{k=1}^n \partial_kf_k=0.$
\textcolor{black}{This is to be understood}
in the generalized function sense that amounts to the relation
%$|\underline{D}|\sum_{k=1}^n R_kf_k=0$ or
\textcolor{black}{$\sum_{k=1}^n R_kf_k=0$} or, equivalently, $\sum_{k=1}^n \xi_k\hat{f_k}(\xi)=0,$ through the inverse Fourier transformation\textcolor{black}{,
again to be understood in the generalized function sense if $p>2$ so that
the Fourier transform is really a distribution. The space of vector valued
divergence free maps  in $L^p(\R^n,\R^n)$
is a closed subset thereof and thus a Banach space in its own right
 that we denote by
$D^p(\R^n)$. Though initially
defined on $\R^n$ only, a divergence free
vector field extends naturally to $\R_{1,\pm}^n$ (resp. $\R^{n+1}_\pm$)
into a $\R^n$ valued map $F$ which is independent of $x_0$ (resp. $x_{n+1}$).
This function needs not be monogenic, but it satisfies
${\rm Sc}\,DF=0$ (resp. ${\rm Sc}\,D_{n+1}F=0$)}.

\section{\textcolor{black}{Hardy-Hodge Decomposition of Para-Vector-Valued functions in $L^p(\mathbb{R}^{n})$}}
\label{sHH}
\begin{theorem}\label{first}
Let \textcolor{black}{$f$ be a para-vector valued function in
$L^p(\mathbb{R}^n,Cl(n,\Phi))$}, $1<p<\infty.$ Then $f$ is uniquely decomposed as
$f=f^++f^-+f^0,$ all in $L^p(\mathbb{R}^n),$ such that $f^\pm$  are para-vector-valued, being the non-tangential boundary limits of some two functions in, respectively, $H^p(\mathbb{R}^{n}_{1,\pm}, \mathbb{R}^n_1),$ and $f^0$ is \textcolor{black}{vector-valued and divergence free}. Moreover, for all $p$ in \textcolor{black}{the indicated} range the decomposition is unique  and, for $p=2,$
\begin{eqnarray}\label{more}
\| f\|^2=\|f^+\|^2+\|f^-\|^2+\|f^0\|^2.\end{eqnarray}
In fact, this decomposition induces a topological direct sum:
\begin{equation}
\label{DS}
L^p(\R^n,\R^n_1)=H^p(\mathbb{R}^{n}_{1,+}, \mathbb{R}^n_1)\oplus
H^p(\mathbb{R}^{n}_{1,-}, \mathbb{R}^n_1)\oplus D^p(\R^n).
\end{equation}
\end{theorem}

\noindent{\bf Proof}
Let $f(\underline{x})=\sum_{k=0}^n f_k(\underline{x})\mathbf{e}_k$ be in $L^p(\mathbb{R}^n),$ where $f_0(\underline{x}),f_1(\underline{x}),\cdots, f_n(\underline{x})$ are scalar-valued and  $\underline{x}=x_1\mathbf{e}_1+\cdots +x_n\mathbf{e}_n\in \mathbf{R}^n.$

By (\ref{Cauchy PI}) \textcolor{black}{and \eqref{Cauchy PIm},
since $H^2=I$}, we have that
\textcolor{black}{$\mathcal{C}^\pm$ is a projection:}
\begin{eqnarray*}
\frac{1}{2}(I+H)f(\underline{x})=[\frac{1}{2}(I+H)]^2f(\underline{x}),\\
\frac{1}{2}(I-H)f(\underline{x})=[\frac{1}{2}(I-H)]^2f(\underline{x}).
\end{eqnarray*}
Then we have
\begin{eqnarray}\label{all}
f(\underline{x})&=&\frac{1}{2}(I+H)f(\underline{x})+\frac{1}{2}(I-H)f(\underline{x})\nonumber \\
&=&[\frac{1}{2}(I+H)]^2f(\underline{x})+[\frac{1}{2}(I-H)]^2f(\underline{x})\nonumber \\
&=&\frac{1}{2}(I+H)\left[{\rm Sc}\{\frac{1}{2}(I+H)f(\underline{x})\}+{\rm Nsc}\{\frac{1}{2}(I+H)f(\underline{x})\}\right]\nonumber \\&+&
\frac{1}{2}(I-H)\left[{\rm Sc}\{\frac{1}{2}(I-H)f(\underline{x})\}+{\rm Nsc}\{\frac{1}{2}(I-H)f(\underline{x})\}\right]\nonumber \\
&=&\frac{1}{2}(I+H)[{\rm Sc}\{\frac{1}{2}(I+H)f(\underline{x})\}]+
\frac{1}{2}(I-H)[{\rm Sc}\{\frac{1}{2}(I-H)f(\underline{x})\}]\nonumber \\
&+&\frac{1}{2}(I+H)[{\rm Nsc}\{\frac{1}{2}(I+H)f(\underline{x})\}]+
\frac{1}{2}(I-H)[{\rm Nsc}\{\frac{1}{2}(I-H)f(\underline{x})\}].
\end{eqnarray}
\textcolor{black}{Consider} the function given by the last line of the above chain of equalities, viz.
\begin{eqnarray}\label{Nsc}\frac{1}{2}(I+H)[{\rm Nsc}\{\frac{1}{2}(I+H)f(\underline{x})\}]+
\frac{1}{2}(I-H)[{\rm Nsc}\{\frac{1}{2}(I-H)f(\underline{x})\}].\end{eqnarray}

It can be computed directly \textcolor{black}{through}
\begin{eqnarray*}
&&\frac{1}{2}(I+H)[{\rm Nsc}\{\frac{1}{2}(I+H)f(\underline{x})\}]\\
&=&\frac{1}{4}(I+H)[\sum_{k=1}^n f_k\mathbf{e}_k+{\rm Nsc}\{Hf\}]\\
&=&\frac{1}{4}\left\{\sum_{k=1}^n f_k\mathbf{e}_k+{\rm Nsc}\{Hf\}+\sum_{k=1}^n H[f_k\mathbf{e}_k]+H[{\rm Nsc}\{Hf\}]\right\},
\end{eqnarray*}
and
\begin{eqnarray*}
&&\frac{1}{2}(I-H)[{\rm Nsc}\{\frac{1}{2}(I-H)f(\underline{x})\}]\\
&=&\frac{1}{4}(I-H)[\sum_{k=1}^n f_k\mathbf{e}_k-{\rm Nsc}\{Hf\}]\\
&=&\frac{1}{4}\left\{\sum_{k=1}^n f_k\mathbf{e}_k-{\rm Nsc}\{Hf\}-\sum_{k=1}^n H[f_k\mathbf{e}_k]]+H[{\rm Nsc}\{Hf\}]\right\}.
\end{eqnarray*}
By adding \textcolor{black}{these relations} together, we have
\textcolor{black}{that}
\begin{eqnarray}
\nonumber
&&\frac{1}{2}(I+H)[{\rm Nsc}\{\frac{1}{2}(I+H)f(\underline{x})\}]+\frac{1}{2}(I-H)[{\rm Nsc}\{\frac{1}{2}(I-H)f(\underline{x})\}]\\
\label{compi}
&=&\frac{1}{2}\left\{\sum_{k=1}^n f_k\mathbf{e}_k+H[{\rm Nsc}\{Hf\}]\right\}.
\end{eqnarray}
\textcolor{black}{As $f$ is para-vector valued,}
relation (\ref{all}) and the fact that $H$ maps scalar valued functins to vector valued functions
together imply that \textcolor{black}{the quantity in \eqref{compi}
is a paravector, therefore}
$H[{\rm Nsc}\{Hf\}]$ is a para-vector. Now we work out its \textcolor{black}{expression}.
Since
\begin{eqnarray*}
Hf=\sum_{k=1}^nR_k(f_k)-\sum_{k=1}^nR_k(f_0){\bf e}_k+ H\wedge\underline{f},
\end{eqnarray*}
where $\underline{f}=\sum_{k=1}^n f_k{\bf e}_k,$ we have \textcolor{black}{that}
\begin{eqnarray*}
{\rm Nsc}\{Hf\}=-\sum_{k=1}^nR_k(f_0){\bf e}_k+ H\wedge\underline{f}.
\end{eqnarray*}
Consequently, \textcolor{black}{since we need only collect terms of degree $0$ and $1$, and because Riesz transforms commute, we obtain:}
\begin{eqnarray*}
H[{\rm Nsc}\{Hf\}]&=& (-\sum_{k=1}^n R_k^2)f_0+\sum_{k=1}^n\left[(\sum_{l\neq k}-R_l^2)f_k +R_k\sum_{l\ne k}R_lf_l\right]{\bf e}_k\\
&=& f_0+\sum_{k=1}^n\left[(\sum_{l\neq k}-R_l^2)f_k +R_k\sum_{l\ne k}R_lf_l\right]{\bf e}_k.
\end{eqnarray*}
Substituting back \textcolor{black}{in}to \textcolor{black}{(\ref{compi})}, we \textcolor{black}{get that}
\begin{eqnarray*}
&&\frac{1}{2}(I+H)[{\rm Nsc}\{\frac{1}{2}(I+H)f(\underline{x})\}]+\frac{1}{2}(I-H)[{\rm Nsc}\{\frac{1}{2}(I-H)f(\underline{x})\}]\\
&=&\frac{1}{2}\left\{f(\underline{x})+\sum_{k=1}^n\left[(\sum_{l\ne k}-R_l^2)f_k +R_k\sum_{l\ne k}R_lf_l\right]{\bf e}_k.\right\}.\end{eqnarray*}

Therefore, \textcolor{black}{by \eqref{all}},

\begin{eqnarray*}
f(\underline{x})
&=&\frac{1}{2}(I+H)\left[{\rm Sc}\{\frac{1}{2}(I+H)f(\underline{x})\}\right]+
\frac{1}{2}(I-H)\left[{\rm Sc}\{\frac{1}{2}(I-H)f(\underline{x})\}\right]\nonumber \\
&&+\frac{1}{2}\left\{f(\underline{x})+\sum_{k=1}^n\left[(\sum_{l\ne k}-R_l^2)f_k +R_k\sum_{l\ne k}R_lf_l\right]{\bf e}_k\right\}.
\end{eqnarray*}

Finally,
\begin{eqnarray}\label{finally}
f(\underline{x})
 &=&(I+H)\left[{\rm Sc}\{\frac{1}{2}(I+H)f(\underline{x})\}\right]+
(I-H)\left[{\rm Sc}\{\frac{1}{2}(I-H)f(\underline{x})\}\right]+\nonumber \\
&&+\sum_{k=1}^n\left[(\sum_{l\ne k}-R_l^2)f_k +R_k\sum_{l\ne k}R_lf_l\right]{\bf e}_k.
\end{eqnarray}
It is apparent that $(I+H)[{\rm Sc}\{\frac{1}{2}(I+H)f(\underline{x})\}]$ and $
(I-H)[{\rm Sc}\{\frac{1}{2}(I-H)f(\underline{x})\}]$ are para-vector-valued, and \textcolor{black}{it follows from Lemma \ref{bvH} that they are
boundary values of functions in $H^p(\R_{1,+}^n,\R_1^n)$ and
$H^p(\R_{1,-}^n,\R_1^n)$ respectively.}
%non-tangential boundary limits of functions in, respectively, the Hardy spaces of the upper- and the lower-half spaces.
Now we show that
\[\sum_{k=1}^n\left[(\sum_{l\ne k}-R_l^2)f_k +R_k\sum_{l\ne k}R_lf_l\right]{\bf e}_k\] is divergence-free.
\textcolor{black}{For this, by} the last remark of the last section, it suffices to show \textcolor{black}{that}
\begin{eqnarray}\label{thus}\sum_{k=1}^nR_k\left[(\sum_{l\ne k}-R_l^2)f_k +R_k\sum_{l\ne k}R_lf_l\right]=0.\end{eqnarray}
The above, \textcolor{black}{however, is} obvious \textcolor{black}{since Riesz transformations commute, and thus
we obtain the} desired decomposition
$f(\underline{x})=f^++f^-+f^0,$
where
\begin{equation}
\label{term1}f^+=(I+H)[{\rm Sc}\{\frac{1}{2}(I+H)f(\underline{x})\}]\in H^p(\mathbb{R}^{n}_{1,+},\R^n_1),
\end{equation}
\begin{equation}
\label{term2}f^-=(I-H)[{\rm Sc}\{\frac{1}{2}(I-H)f(\underline{x})\}]\in H^p(\mathbb{R}^{n}_{1,-},\R^n_1),
\end{equation}
and
\begin{equation}
\label{term3}
f^0=\sum_{k=1}^n\left[(\sum_{l\ne k}-R_l^2)f_k +R_k\sum_{l\ne k}R_lf_l\right]{\bf e}_k
\end{equation}
is divergence free. \textcolor{black}{Next we
prove}  uniqueness. \textcolor{black}{This} is equivalent \textcolor{black}{to} showing that if we have a decomposition of the zero function\textcolor{black}{:} $0=f^++f^-+f^0,$ then \textcolor{black}{it must be that} $f^+=f^-=f^0=0.$
Indeed, in \textcolor{black}{that} case we \textcolor{black}{ may write}
\[ 0=(I+H)(f^++f^-+f^0)=2f^++(I+H)f^0=2f^++f^0+Hf^0,\]
\textcolor{black}{where we used that $f^+=Hf^+$ by \eqref{Cauchy PI} since
$\mathcal{C}^+f^+=f^+$, and also that
$Hf^-=-f^-$ by \eqref{Cauchy PIm} since $\mathcal{C}^-f^-=f^-$.}
Note that since $f^0$ is divergence free, the scalar part of $Hf^0$ is zero, and thus only the $2$-form part of $Hf^0$ is possibly non-zero. However, the
last \textcolor{black}{equality}  shows that the $2$-form part also has to be zero, because \textcolor{black}{all the other terms are} para-vectors. We thus conclude
\textcolor{black}{ that} $2f^++f^0=0.$ The same reasoning gives
\textcolor{black}{us}
$2f^-+f^0=0.$ \textcolor{black}{These} together \textcolor{black}{yield} $f^+=f^-.$ By applying $I+H$ \textcolor{black}{to} both sides we \textcolor{black}{get}
 $f^+=0$ \textcolor{black}{hence also} $f^-=0$, \textcolor{black}{and
consequently}
 $f^0=0.$ \textcolor{black}{This establishes uniqueness and shows that
\eqref{DS} holds as a direct sum.
In addition,
since $f^\pm$ and $f^0$ are continuous functions of $f$ in $L^p(\R^n, Cl(n,\Phi))$ by
\eqref{term1}, \eqref{term2},
\eqref{term3} and the $L^p$ continuity of the Riesz transformations,
we see that the projections in \eqref{DS} are continuous, hence the sum is topological by the open mapping theorem.}
Finally, \textcolor{black}{when $p=2$,}  we show the \textcolor{black}{Pythagora} type relation
\begin{eqnarray}\label{orthogonal pro}\|f\|^2
=\|f^+\|^2+\|f^-\|^2+\|f^0\|^2.\end{eqnarray}
 First, since $f$ is para-vector-valued, \textcolor{black}{we obviously have that}
\[ \| f\|^2=\int_{\mathbf{R}^2}f \overline{f} d\underline{x}.\]
\textcolor{black}{Hence,} to prove (\ref{orthogonal pro}), it suffices  to
\textcolor{black}{establish} the following \textcolor{black}{orthogonality}
relations:
\begin{eqnarray}\label{equality 1}\int_{\mathbf{R}^n}f^+\overline{f^-}d\underline{x}=\label{equality 2}\int_{\mathbf{R}^n}f^-\overline{f^+}d\underline{x}
=0,\end{eqnarray}
and
\begin{eqnarray}\label{equality 3}\int_{\mathbf{R}^n}(f^+\overline{f^0}+f^0\overline{f^+})d\underline{x}=
\int_{{\bf R}^n}(f^-\overline{f^0}+f^0\overline{f^-})d\underline{x}
=0.\end{eqnarray}
\textcolor{black}{Let us show} (\ref{equality 1}).
\textcolor{black}{Recall from Lemma \ref{bvH}} that if $g\in H^2(\mathbb{R}^{n}_{1,\pm}, \mathbb{R}^n_1),$ then $g=(I\pm H)g_0$ and  \textcolor{black}{consequently, taking Fourier transforms (the Fourier transform of a vector valued
function is computed componentwise), we get that}
$\hat{g}=\textcolor{black}{2}\chi_\pm\hat{g_{\textcolor{black}{0}}},$ where $\chi_\pm $ are \textcolor{black}{multipliers for} the Hardy space projections\textcolor{black}{:} $\chi_\pm (\underline{\xi})=\frac{1}{2}(1\pm i\frac{\underline{\xi}}{|\underline{\xi}|}),$
\textcolor{black}{that satisfy} $\chi_\pm^2=\chi_\pm$ \textcolor{black}{and} $\chi_++\chi_-=1$
\textcolor{black}{as well as} $\chi_+\chi_-=\chi_-\chi_+=0$. \textcolor{black}{Here we used the
expression for the multiplier of $R_k$ in the Fourier domain, see}
\cite{St,LMcQ}.   \textcolor{black}{Applying these} remarks to $g=f^\pm,$
\textcolor{black}{and using} Parseval's Theorem, we have \textcolor{black}{that}
\begin{eqnarray*}&&\int_{\mathbf{R}^2}f^+\overline{f^-}d\underline{x}\\
&=&\int_{\mathbf{R}^2}(f^+)^\wedge(\underline{\xi})\overline{(f^-)^\wedge(\underline{\xi})}d\underline{\xi}
\\
%&=&\int_{\mathbf{R}^2}2\chi_+(\underline{\xi})[{\rm Sc}\{\frac{1}{2}(I+H)f(\cdot)\}]^\wedge(\underline{\xi})\overline{2\chi_-(\underline{\xi})[{\rm Sc}\{\frac{1}{2}(I-H)f(\cdot)\}]^\wedge(\underline{\xi})}d\underline{\xi}\\
&=&\textcolor{black}{4\int_{\mathbf{R}^2}\chi_+(\underline{\xi})\chi_-(\underline{\xi})
%[{\rm Sc}\{\frac{1}{2}(I+H)f(\cdot)\}]^\wedge(\underline{\xi})\overline{[{\rm Sc}\{\frac{1}{2}(I-H)f(\cdot)\}]^\wedge(\underline{\xi})}
\hat{f_0^+}\overline{\hat{f_0^-}}\,d\underline{\xi}}\\
&=&0,
\end{eqnarray*}
where we used the relation $\chi_+(\underline{\xi})\chi_-(\underline{\xi})=0$ for all $\xi.$ The proof of the second equality relation in (\ref{equality 2}) is similar.
\textcolor{black}{Now, let us show that} equality (\ref{equality 3})
\textcolor{black}{holds}.
Indeed,
\begin{eqnarray*}
\int_{\mathbf{R}^n}(f^+\overline{f^0}+f^0\overline{f^+})d\underline{x}&=&
2{\rm Sc}\left\{\int_{\mathbb{R}^n} f^+\overline{f^0}d\underline{x}\right\}\\
&=& {\rm Sc}\left\{\int_{\mathbb{R}^n}\left(1+i\frac{\underline{\xi}}{|\underline{\xi}|}\right)
{\widehat{f_0}}(\xi)\overline{{\widehat{f^0}}}(\underline{\xi})
d\underline{\xi}\right\}\\
&=&  {\rm Sc}\left\{\int_{\mathbb{R}^n}{\widehat{f_0}}(\xi)\overline{\frac{i\underline{\xi}}
{|\underline{\xi}|}{{\widehat{f^0}}}(\underline{\xi})}
d\underline{\xi}\right\}\\
&=&0,
\end{eqnarray*}
where the last equality used the relation \textcolor{black}{
$(i\underline{\xi}/
|\underline{\xi}|){\widehat{f^0}}(\underline{\xi})=0,$ the latter being
 a consequence of the fact that} $f^0$ is divergence free.
The proof is complete.

\noindent{\bf Remark} \textcolor{black}{If, alternatively, we use the scalar product}
\[ \langle f,g\rangle = {\rm Sc} \int_{\mathbb{R}^n} f \overline{g}d\underline{x},\]
then indeed the decomposition $f=f^++f^-+f^0$ is \textcolor{black}{orthogonal}.
\textcolor{black}{Moreover, since we observed after the proof of Lemma \ref{bvH}
that $H^2(\R^n_{1,\pm},\R_1^n)\cap H^p(\R^n_{1,\pm},\R_1^n)$  is dense in
$H^p(\R^n_{1,\pm},\R_1^n)$, we deduce that \eqref{equality 1} holds as
soon as $f^+\in H^p(\R^n_{1,+},\R_1^n)$ and
$f^-\in H^{p'}(\R^n_{1,-},\R_1^n)$  with $1/p+1/p'=1$:
\begin{equation}
\label{equality1p}
. \int_{\mathbf{R}^n}f^+\overline{f^-}d\underline{x}=\int_{\mathbf{R}^n}f^-\overline{f^+}d\underline{x}
=0,\qquad f^+\in H^p(\R^n_{1,+},\R_1^n),\quad
f^-\in H^{p'}(\R^n_{1,-},\R_1^n)
%,\quad 1/p+1/p'=1
.
\end{equation}
Likewise, \eqref{equality 3} generalizes to
\begin{equation}
\label{equality2p}
\int_{\mathbf{R}^n}(f^+\overline{f^0}+f^0\overline{f^+})d\underline{x}=
\int_{{\bf R}^n}(f^-\overline{f^0}+f^0\overline{f^-})d\underline{x}
=0, \quad f^\pm\in H^p(\R^n_{1,\pm},\R_1^n),\quad
f^0\in D^{p'}(\R^n).
\end{equation}}

\textcolor{black}{A few comments are in order:
\begin{itemize}
\item If \eqref{DS} gets projected onto the last $n$ components, and since
$(\pm R_1 h,\cdots,\pm R_n h)$ is a gradient vector field on $\R^n$
({\it i.e.} the gradient of the trace  of a solution to the Neumann
problem on $\R^{n+1}_\pm$ with inner normal derivative $h$ a.e. on $\R^n$)
we recover the classical Helmoltz-Hodge decomposition of vector fields from
$L^p(\R^n,\R^n)$ into the sum of a rotational-free and a divergence-free
vector field \cite{im2}.
\item Decomposition \eqref{DS} generalizes to higher dimensions
the standard decomposition of a complex valued function in $L^p(\R)$
into the sum of
a function belonging to the holomorphic Hardy space $H^p(\R^2_+)$ and
a function belonging to the holomorphic Hardy space $H^p(\R^2_-)$. The difference in dimension bigger than 1 is that a divergence free term must be added,
for in this case  not every vector field is a gradient. Note, since $R_k$ and
the divergence operator preserve realness, that Theorem \ref{first} carries
over to
Clifford valued maps and Clifford Hardy spaces with complex coefficients.
\end{itemize}}
\section{Variations}
\label{var}
Next we consider the homogeneous case on $\mathbb{R}^n.$ \textcolor{black}{We
regard $\R^n$ as being the subspace $\R^n\times\{0\}$ of
$\R^{n+1}$. When considering Clifford 1-forms as
Euclidean vectors, it means that 1-forms in $Cl(n,\R)$ get identified with
1-forms in $Cl(n+1,\R)$ whose coefficient of ${\bf e}_{n+1}$ is zero.
Note that $Cl(n,\R)$ can be viewed as the subalgebra of $Cl(n+1,\R)$
generated by ${\bf e}_1,\cdots,{\bf e}_n$.
}

\begin{theorem}\label{second}
Let $f\in L^p(\mathbb{R}^n\textcolor{black}{,Cl(n+1,\R)})$ be an $(n+1)$-vector-valued function, $1<p<\infty.$ Then $f$ is uniquely decomposed as
$f=f^++f^-+f^0,$ all in $L^p(\mathbb{R}^n),$ such that $f^\pm$ are
\textcolor{black}{the} non-tangential boundary limits of some two functions in, respectively, $H^p(\mathbb{R}^{n+1}_\pm,\mathbb{R}^{n+1}),$ \textcolor{black}{while} $f^0$ is \textcolor{black}{vector-valued in $Cl(n,\R)$ } and divergence free. For all $p$ in the \textcolor{black}{indicated range,} the decomposition is
unique \textcolor{black}{and topological}.
Moreover, for $p=2,$ $$\| f\|^2=\| f^+\|^2 + \|f^-\|^2 + \| f^0\|^2.$$
\end{theorem}

 By factorizing out ${\bf e}_{n+1}$ one can, in particular, reduce the proof of Theorem \ref{second} to that of Theorem \ref{first}. \textcolor{black}{More}
precisely, \textcolor{black}{noting that ${\bf e}_{n+1}^{-1}=-{\bf e}_{n+1}$,}
we use the relation

\[ \sum_{k=1}^{n+1} f_k{\bf e}_k=[\sum_{k=1}^{n+1}f_k{\bf e}_k{\bf e}_{n+1}^{-1}]{\bf e}_{n+1}.\]
For the Cauchy kernels in the two settings, one has
\[ \frac{\overline{\sum_{k=1}^ny_k\mathbf{e}_k-(\sum_{k=1}^{n+1}x_k\mathbf{e}_k)}}
{|\sum_{k=1}^n{y_k\mathbf{e}_k-(\sum_{k=1}^{n+1}x_k\mathbf{e}_k)}|^{n+1}}=-{\bf e}_{n+1} \frac{\overline{\sum_{k=1}^n y_k\mathbf{e}_k
{\bf e}_{n+1}^{-1}-(\sum_{k=1}^{n}x_k\mathbf{e}_k{\bf e}_{n+1}^{-1}+x_{n+1})}}
{|{\sum_{k=1}^n y_k\mathbf{e}_k{\bf e}_{n+1}^{-1}-(\sum_{k=1}^n x_k\mathbf{e}_k{\bf e}_{n+1}^{-1}+x_{n+1})}|^{n+1}}.\]
\textcolor{black}{This corresponds to} the relation between the two Dirac operators:
\[ \sum_{k=1}^{n+1} \partial_k{\bf e}_k=(\sum_{k=1}^{n+1} \partial_k{\bf e}_k{\bf e}_{n+1}^{-1}){\bf e}_{n+1}.\]
\textcolor{black}{Letting} $\tilde{\bf e}_k={\bf e}_k{\bf e}_{n+1}^{-1}, \ k=1,...,n,$ one reduces the proof of Theorem \ref{second} to that of Theorem \ref{first} (\textcolor{black}{compare} \cite{Q}).

\textcolor{black}{Theorem} \ref{second} can alternatively \textcolor{black}{be re}written \textcolor{black}{without mentioning Clifford analysis} as

\begin{theorem}\label{second plus}
Every vector field $f\in L^p(\mathbb{R}^n, \mathbb{R}^{n+1}),$  $1<p<\infty,$ may be uniquely decomposed \textcolor{black}{as} $f=f^++f^-+f^0,$ \textcolor{black}{ where $f^\pm\in L^p(\mathbb{R}^n, \mathbb{R}^{n+1})$} are, respectively, the non-tangential boundary limits of some harmonic gradients \textcolor{black}{on}
${\mathbb{R}}^{n+1}_\pm,$ \textcolor{black}{which} satisfy  (\ref{gradiants}),
\textcolor{black}{while  $f^0\in L^p(\R^n,\R^n)$ is} divergence free.
\textcolor{black}{The decomposition is topoligical and} for $p=2$ there holds
\[ \| f\|^2=\| f^+\|^2 + \|f^-\|^2 + \| f^0\|^2.\]
\end{theorem}
\textcolor{black}{From the analog of Lemma \ref{bvH} for
homogeneous Hardy spaces (see the discussion after the proof of that lemma),
we know that  boundary limits of harmonic gradients on
${\mathbb{R}}^{n+1}_\pm$  satisfying  (\ref{gradiants}) are those
members of $L^p(\R^n,\R^{n+1})$ of the type
$(\pm R_1 h,\cdots,\pm R_n h, h)$ with  $h\in L^p(\R^n, \R),$ where all the $\lq\lq$+" signs and, respectively, all the $\lq\lq$-" signs are taken.
In such form, Theorem \ref{second plus}
was proven in \cite{BHLSW} when $n=2$,
also for more general function spaces.}

\textcolor{black}{It is worth contrasting  Theorem \ref{first} and Theorem \ref{second} with their quaternionic counterparts. The space $\mathbb{H}$ of
real quaternions consists of numbers}
$q=q_0+\underline{q},$ \textcolor{black}{with}
$\underline{q}=q_1{\bf e}_1+q_2{\bf e}_2+q_3{\bf e}_3$ \textcolor{black}{where
$q_j\in\R$ and } ${\bf e}_1, \ {\bf e}_2$
\textcolor{black}{are as before, but, additionally},
${\bf e}_3={\bf e}_1{\bf e}_2.$ We identify the linear space consisting of all $\underline{q}$ with the space $\mathbb{R}^3$, \textcolor{black}{and we put
$\mathbb{H}_\pm$ for those quaternions with, respectively, $\pm q_0>0$. We say that
a quaternionic valued function $f$  is left quaternionic if
$(D_0+D_3)f=0$, see \eqref{defop},
but this time the relation
${\bf e}_1{\bf e}_2={\bf e}_3$ is taken into account. The definitions of
quaternionic Hardy spaces $H^p(\mathbb{H}_\pm, \mathbb{H})$ as spaces of
left quaternionic functions in $\mathbb{H}_\pm$ meeting the analogs of
\eqref{firste} now run parallel to those for inhomogeneous Hardy spaces.
}
\begin{theorem}\label{third}
Let $f\in L^p(\mathbb{R}^3, \mathbb{H}), \ 1<p<\infty.$ Then $f$ \textcolor{black}{is uniquely decomposed} as
$f=f^++f^-$ such that \textcolor{black}{$f^\pm\in L^p(\mathbb{R}^3, \mathbb{H})$} are non-tangential boundary limit functions of some two functions in, respectively, $H^p(\mathbb{H}_\pm, \mathbb{H}).$ Moreover, for $p=2,$ $f^+$ and $f^-$ are orthogonal\textcolor{black}{:}
\[ \| f\|^2 = \|f^+ \|^2 + \| f^- \|^2.\]
\end{theorem}

\noindent{Proof\textcolor{black}{.}} This is an immediate consequence of the corresponding Plemelj formula, for functions obtained from the Cauchy formula are all
\textcolor{black}{quaternionic}-valued.\\

\noindent{\bf Remark\textcolor{black}{:}} \textcolor{black}{the reason why the divergence free term} $f^0$ \textcolor{black}{can  be omitted} in Theorem \ref{third} is the close\textcolor{black}{d}ness \textcolor{black}{of multiplication} in the 
\textcolor{black}{quaternionic} field. \textcolor{black}{At the same time, the interpretation of $f^\pm$ as traces of harmonic gradients is lost}.
To see the difference \textcolor{black}{in} the analysis \textcolor{black}{of} $H^p(\mathbb{R}^{3}_{1,\pm}, \mathbb{R}^3_1 )$ and $H^p(\mathbb{H}_\pm, \mathbb{H})$\textcolor{black}{, recall that if} $f(x_0+\underline{x})=f_0(x_0+\underline{x}){\bf e}_0+f_1(x_0+\underline{x}){\bf e}_1+f_2(x_0+\underline{x}){\bf e}_2+f_3(x_0+\underline{x}){\bf e}_3$ \textcolor{black}{lies in} $H^p(\mathbb{R}^{3}_{1,+}, \mathbb{R}^3_1),$ \textcolor{black}{ then}
\begin{equation}\label{CHS for clifford}
\begin{cases}\frac{\partial f_0}{\partial x_0}=\frac{\partial f_1}{\partial x_1}+\frac{\partial f_2}{\partial x_2}+\frac{\partial f_3}{\partial x_3} ,\\ \frac{\partial f_0}{\partial x_i}=-\frac{\partial f_i}{\partial x_0},\ \ i=1,2,3,\\
\frac{\partial f_i}{\partial x_j}=\frac{\partial f_j}{\partial x_i},\ \ i\ne 0, j\ne 0, i\neq j,\end{cases}\end{equation}
while
$f(q_0+\underline{q})=f_0(q_0+\underline{q}){\bf e}_0+f_1(q_0+\underline{q}){\bf e}_1+f_2(q_0+\underline{q}){\bf e}_2+f_3(q_0+\underline{q}){\bf e}_3\in H^p(\mathbb{H}_+, \mathbb{H})$ will imply
\begin{equation}\label{CHS for quoter}
\begin{cases}\frac{\partial f_0}{\partial x_0}=\frac{\partial f_1}{\partial x_1}+\frac{\partial f_2}{\partial x_2}+\frac{\partial f_3}{\partial x_3} ,\\\frac{\partial f_1}{\partial x_2}=\frac{\partial f_3}{\partial x_0}+\frac{\partial f_0}{\partial x_3}+\frac{\partial f_2}{\partial x_1}\\
\frac{\partial f_2}{\partial x_3}=\frac{\partial f_1}{\partial x_0}+\frac{\partial f_0}{\partial x_1}+\frac{\partial f_3}{\partial x_2}\\
\frac{\partial f_3}{\partial x_1}=\frac{\partial f_2}{\partial x_0}+\frac{\partial f_0}{\partial x_2}+\frac{\partial f_1}{\partial x_3}\end{cases}.
\end{equation}

The system of equations (\ref{CHS for clifford}) implies the system of equations (\ref{CHS for quoter}), \textcolor{black}{but not conversely}.

\section{\textcolor{black}{Uniqueness issues for potentials in divergence form}}
\label{sis}
\textcolor{black}{In \cite{BHLSW}, the Hardy-Hodge decomposition was introduced when $n=2$ to
characterize silent magnetizations with support in $\R^2$. This issue can be
recast as that of describing vanishing potentials in divergence form.
Recall that the Newton potential of a distribution $\varphi$
on $\R^{n+1}$ is the convolution of $\varphi$ with $1/(\omega_n|x|^{n-1})$,
wherever it exists.
% \begin{equation}
% \label{NP}
% P_\varphi(x)=\frac{1}{\omega_n}\int \frac{1}{|x-y|^{n-1}} \varphi(y)
% \end{equation}
% where $\varphi$ is assumed smooth enough that \eqref{NP} makes sense at least for $x\notin{\rm supp}\,\varphi$, where ${\rm supp}\, \varphi$ indicates the support of $\varphi$.
The potential is said to be in divergence form if
$\varphi$ can be taken to be the divergence of some $\R^{n+1}$ valued
distribution $\psi$:
\begin{equation}
\label{NPD}
P_{{\rm div}\,\psi}(x)=\frac{1}{\omega_n}\int \frac{1}{|x-y|^{n-1}} {\rm div}\,\psi(y)=-\frac{1}{\omega_n}\int \nabla_y\left(\frac{1}{|x-y|^{n-1}}\right)
\cdot \psi(y),
\end{equation}
where $\nabla_y$ indicates the gradient with respect to the variable $y,$ where the dot indicates Euclidean scalar product.
When $n=2$, in the quasi-static approximation to Maxwell's equations,
\eqref{NPD} formally expresses the magnetic potential of the
magnetization $\psi$ (cf. \cite[Section 5.9.C]{Jackson}).
Those $x$ (if any) for which this expression makes good sense depend of
course on
$\psi$. We shall be concerned with the case where $\psi$ is supported
on a hyperplane $\mathcal{P}$ and has $L^p$ density there.
Specifically, if we write
\begin{equation}
\label{defhyp}
\mathcal{P}=\{x\in \R^{n+1},\ x\cdot u=a\}
\end{equation}
for some $u\in\R^{n+1}$ and $a\in\R$,
 it means that
$\psi=f\otimes\delta_0(x\cdot u-a)$ where
$f=(f_0,\cdots,f_{n})\in L^p(\mathcal{P},\R^{n+1})$  and $\delta_0$ indicates
the Dirac mass at $0$ .
Then, \eqref{NPD} becomes
\begin{equation}\label{tpphip}
P_{{\rm div}\,\psi}(x)=
\frac{n-1}{\omega_n}\int_{\mathcal{P}}
\frac{f(y)\cdot (x-y)}{|x-y|^{n+1}}
\, dy,\end{equation}
which is well defined for all $x\notin\mathcal{P}$, more generally
for all $x$ not in the support of $f$.}

\textcolor{black}{Let
\[
\mathcal{H}_\pm=\{x\in \R^{n+1},\ \pm(x\cdot u-a)>0\}\]
denote the two half spaces whose union is $\R^{n+1}\setminus\mathcal{P}$,
the complement of $\mathcal{P}$.
The question that we raise is:
\begin{enumerate}
\item[] \emph{for which $f$ does it happen that
$P_{{\rm div}\,\psi}(x)=0$ for all $x\in\mathcal{H}_\pm$?}
\end{enumerate}
From the physical viewpoint, it amounts when $n=1$
to describe those magnetizations
with $L^p$ density supported on a plane  which are silent from one side of
that plane, meaning that they generate no magnetic field in the corresponnding half-space.
These cannot be detected by measuring devices and account a good deal for
the ill-posedness of inverse magnetization problems \cite{Parker}.
The result below gives an answer to the question in terms of the Hardy-Hodge decomposition, thereby generalizing to higher dimension results from
\cite{BHLSW} for $n=2$.}

\textcolor{black}{It will be convenient to define the Clifford Hardy spaces
$H^p(\mathcal{H}_\pm,\R^n_1)$ consisting of para-vector valued monogenic
functions $g$ in $\mathcal{H}_\pm$ meeting the condition:
\begin{equation}
\label{charHP}
\sup_{\pm b>0} \int_{\mathcal{P}_b}|g|^p\,dm<\infty,\qquad \mathcal{P}_b\triangleq \{x\in \R^{n+1},\ \pm(x\cdot u-a)=b\},
\end{equation}
where $dm$ indicates the differential of Lebesgue measure.
Just like in the case of inhomogeneous Hardy spaces $H^p(\R^n_{1,\pm},\R^n_1)$,
condition \eqref{charHP} may be replaced by the
$L^p(\mathcal{P},\R)$-boundedness of the nontangential maximal function,
computed this time over cones with vertex on $\mathcal{P}$ and
axis parallel to $u$. Functions in $H^p(\mathcal{H}_\pm,\R^n_1)$
have nontangential limits in
$L^p(\mathcal{P},\R_1^n)$  of which they are the Cauchy Clifford
integral, and they can be identified with their nontangential limit.
In fact, if
$\mathfrak{R}$ is any orientation preserving affine isometry mapping
$\mathcal{P}$ to $\{0\}\times\R^n$, we have that $f$ belongs to
$H^p(\mathcal{H}_\pm,\R^n_1)$ if and only if $f\circ\mathfrak{R}$ belongs to
$H^p(\R^n_{1,\pm}, \R^n_1)$. As an extra piece of notation, we put $u^\perp$
for the vector space orthogonal to $u$ in $\R^{n+1}$, which is the linear
hyperplane parallel to $\mathcal{P}$.
\begin{theorem}
Let $\mathcal{P}\subset\R^{n+1}$ be a hyperplane \textcolor{black}{defined by
\eqref{defhyp}} and $\psi\in L^p(\mathcal{P},\R^{n+1})$ with $1<p<\infty$. Then the potential $P_{{\rm div}\psi}$ vanishes on
$\mathcal{H}_\pm$ if and only if $\psi$ is the sum  of a member of
$H^p(\mathcal{H}\pm,\R^n_1)$ and of a divergence free \textcolor{black}{vector field tangent to $\mathcal{P}$} in
$L^p(\mathcal{P}, u^\perp)$. The potential $P_{{\rm div}\psi}$ vanishes on
$\R^{n+1}\setminus\mathcal{P}$ (that is, on both
$\mathcal{H}_+$ and $\mathcal{H}_-$) if and only if it is
a divergence free function in
$L^p(\mathcal{P}, u^\perp)$.
\label{SS}
\end{theorem}
\noindent{Proof.}
Because the statement is invariant under orientation preserving
affine isometries, we may assume that $\mathcal{P}=\{0\}\times\R^n$ so that
$\mathcal{H}_\pm=\R^n_{1,\pm}$. Let us single out the components of
$\psi$ as $\psi_0,\cdots,\psi_{n}$ and identify $\psi$ with the para-vector
valued  function $\psi=\psi_0{\bf e}_0+\psi_1{\bf e}_1+\cdots+\psi_n{\bf e}_n$.
Set $\psi=\psi^++\psi^-+\psi^0$ for the Hardy-Hodge decomposition from Theorem
\ref{first}. For $x\in\R^n_{1,\pm}$, it is
easily checked that $y\mapsto x-y/|x-y|^{n+1}$ lies in
$H^q(\R^n_{1,\pm},\R^n_1)$. Thus, it follows from \eqref{tpphip},
\eqref{equality1p} and \eqref{equality2p} that
$P_{{\rm div}(\psi^\mp+\psi_0)}\equiv0$ on $\R^n_{1,\pm}$. Therefore,
the assumption that  $P_{{\rm div}\psi}\equiv0$ on $\R^n_{1,\pm}$ reduces to
$P_{{\rm div}\psi^\mp}\equiv0$ on $\R^n_{1,\pm}$.
Now, comparing \eqref{tpphip} and \eqref{defCauchyp}, we find this is
equivalent to
\[
{\rm Sc}\left\{ C^\pm \psi^\mp\right\}(\bar{x})=0,\qquad x\in\R^n_{1,\pm}.
\]
which amounts to
\begin{equation}\label{tpphiC}
{\rm Sc}\left\{ C^\mp \psi^\mp\right\}(\xi)=0,\qquad \xi\in\R^n_{1,\mp}.
\end{equation}
Since $\psi^\mp\in H^p(\R^n_{1,\mp},\R^n_1)$, we have
by the Cauchy Clifford formula that $C^\mp\psi^\mp(\xi)=\psi^\mp(\xi)$,
therefore \eqref{tpphip} means that ${\rm Sc}\{\psi^\mp\}$ vanishes on
$\R^n_{1,\mp}$ and so does its nontangential limit on $\R^n$. But we
know from Lemma \ref{bvH} that the $L^p(\R^n)$-norm of the nontangential
limit of the scalar part is an equivalent norm on $H^p(\R^n_{1,\mp},\R^n_1)$,
hence
$\psi^\mp=0$. This proves the first assertion of the theorem. To establish the
second assertion, observe from what precedes that if $P_{{\rm div}\,\psi}=0$
both in $\R^n_{1,+}$ and $\R^n_{1,-}$, then $\psi^\pm=0,$ and thus
$\psi=\psi^0$ is vector-valued and divergence free. The proof is complete.}

\section*{Acknowledgements}
This work was supported by Macao Science and Technology Development Fund, MSAR. Ref. 045/2015/A2,The work was supported by
Multi-Year Research Grant (MYRG) MYRG116(Y1-L3)-FST13-QT. Macao Government FDCT 098/2012.

\end{document}